\documentclass[10pt,a4paper]{article}
\usepackage[utf8]{inputenc}
\usepackage{amsmath}
\usepackage{amsfonts}
\usepackage{amssymb}
\usepackage{xcolor}

\usepackage{amsthm}
\usepackage{bbm}
\usepackage{graphicx}
\usepackage{float}
\usepackage{mathtools}

\newtheorem{theorem}{Theorem}[section]

\newtheorem{proposition}[theorem]{Proposition}

\newtheorem{conjecture}[theorem]{Conjecture}
\theoremstyle{definition}

\newtheorem{definition}[theorem]{Definition}
\newtheorem{remark}[theorem]{Remark}

\newcommand{\R}{{\mathbb{R}}}
\newcommand{\C}{{\mathbb{C}}}

\newcommand{\ehzcap}{c_{_{\rm EHZ}}}

\newcommand{\dg}{\dot{\gamma}}

\begin{document}

\title{A Counterexample to Viterbo's Conjecture}
\author{Pazit Haim-Kislev, Yaron Ostrover}
\maketitle

\begin{abstract}
We present a counterexample to Viterbo's volume-capacity conjecture.    
This implies, in particular, that in contrast with a well-known conjecture, symplectic capacities do not coincide on the class of convex domains in the classical phase space.
\end{abstract}

\section{Introduction and Results}

A major distinction between the categories of volume-preserving transformations and symplectic transformations was first demonstrated by Gromov in~\cite{gromov}. In this work Gromov introduced the  theory of pseudoholomorphic curves and, in particular, proved the renowned ``non-squeezing" theorem, which provides the first example of obstructions to symplectic embeddings distinct from volume considerations. This discovery later led to the development of global symplectic invariants, known  as {\it symplectic capacities}. 
The axiomatic definition of symplectic capacities 
was first given by Ekeland and Hofer in~\cite{hofer-ekeland}, and since then many examples were devised, arising
 from various themes and techniques in symplectic geometry and Hamiltonian dynamics. These invariants play a central role in symplectic geometry, particularly in establishing rigidity results and identifying obstructions to symplectic embeddings.
Despite their 
importance, many fundamental questions regarding symplectic capacities remain unresolved due to the fact that these capacities are notoriously  difficult to compute.

\medskip 

One of the long-standing fundamental conjectures regarding symplectic capacities is  
Viterbo's volume-capacity conjecture \cite{viterbo2000} (Conjecture \ref{viterbo_conj} below), an isoperimetric-type question concerning convex domains in the classical phase space. Roughly speaking, it  states that among convex domains with a given volume, the Euclidean ball has the maximal symplectic size. This conjecture inspired extensive research, and was proved in many special cases 
(a partial list of relevant references includes~\cite{abbondandoloBenedetti,  abbondandolo-benedetti-edtmair,abhs, ArtsteinAvidanOstroverMilman, balitsky,Ch-H,edtmair, gutt-hutchings-ramos,hermann, karasev-sharipova, rudolf}). 
In particular, in \cite{capacity_mahler}, it was proved that a special case of Viterbo's conjecture is equivalent to the famous Mahler's conjecture, an 85-year-old open problem in convex geometry concerning the volume of a centrally symmetric convex body and that of its dual. 
The aforementioned  works, along with several other results, provided substantial evidence supporting Viterbo's conjecture. 
However, in this paper we provide  a counterexample
based on Minkowski billiard dynamics. We  turn now to a precise statement of the main result.

\medskip

We begin by recalling the notion of symplectic capacity. Various definitions of symplectic capacities exist in the literature, typically defined for the class of symplectic manifolds, or specific subclasses (see, e.g.,~\cite{capacity_survey_1}).
For the purposes of this work, the following definition of what is commonly referred to as a normalized symplectic capacity will suffice.
Consider the classical phase space $\R^{2n} \simeq \C^n$ equipped with canonical  coordinates $(q,p)$, and with the standard symplectic form $\omega_{\rm std} = \sum dp_i \wedge dq_i$.
\begin{definition}
A symplectic capacity is a map $c$ that assigns a non-negative number 
to each subset $U \subset \R^{2n}$, and satisfies the following:
\begin{enumerate}
\item (Monotonicity) $c(U) \leq c(V)$ if there is a symplectic embedding $U \hookrightarrow V$.
\item (Conformality) $c( \alpha U) = \alpha^2 c(U)$, for every $\alpha \in \R_{+}$.
\item (Normalization) $c(B^{2n}(1)) = c(Z^{2n}(1)) = \pi$.
\end{enumerate}
\end{definition}
Here $B^{2n}(r) \subset \C^n$ is the Euclidean ball of radius $r$, and $Z^{2n}(r)$ is the cylinder $B^2(r) \times \C^{n-1}$. 
Two immediate examples of symplectic capacities which arise naturally from Gromov's non-squeezing theorem~\cite{gromov} are the Gromov width 
$$c_{\mathrm{G}}(U) := \sup \{ \pi r^2 : \exists \,  B^{2n}(r) \overset{\mathrm s}{\hookrightarrow } U  \},$$
where $\overset{\mathrm s}{\hookrightarrow }$ stands for a symplectic embedding, and the cylindrical capacity
$$ c_{\mathrm{Z}}(U) = \inf \{ \pi r^2 : \exists \, U  \overset{\mathrm s}{\hookrightarrow } Z^{2n}(r)  \}.$$
Since~\cite{gromov}, numerous additional symplectic capacities (not necessarily satisfying the Normalization property) have been constructed. 
Examples include the Hofer--Zehnder capacity \cite{hofer-zehnder}, the Ekeland--Hofer capacities \cite{hofer-ekeland}, the symplectic homology capacity defined by Hermann \cite{hermann2000} based on ideas by Viterbo \cite{viterbo-sh} (cf. \cite{FHW94}), the Gutt--Hutchings capacities \cite{Gu-Ha} (which were very recently shown~\cite{GuttRamos} to coincide with the Ekeland--Hofer capacities for star-shaped domains), the embedded contact homology capacities \cite{ech} defined by Hutchings in dimension 4, and capacities defined by Siegel using rational symplectic field theory \cite{siegel}.

\medskip

Due to their monotonicity and conformality properties, symplectic capacities are continuous with respect to the Hausdorff topology on the class of smooth convex domains in \(\mathbb{R}^{2n}\). 
Throughout this text, the term convex domain refers to a bounded convex set with nonempty interior.
Consequently, each such capacity extends uniquely to the broader class of convex domains, regardless of smoothness. Moreover, 
for the class of convex domains 
it is known that, among others, the Hofer--Zehnder capacity, the first Ekeland--Hofer capacity, and the symplectic homology capacity 
coincide, and are equal to the minimal action among closed characteristics (or closed Reeb orbits) on the boundary\footnote{Without the smoothness assumption, one needs to consider ``generalized" closed characteristics on the boundary (see, e.g., \textsection 2.3 in~\cite{AA-O}).
} (see, e.g., \cite{A-E-K, hofer-zehnder, hofer-ekeland, viterbo_cap, irie, abbondandolo-kang, Gu-Ha, ginzburg-shon, gutt-hutchings-ramos}). In what follows we refer to this quantity as the Ekeland--Hofer--Zehnder capacity and denote it by $\ehzcap$. As mentioned above, an influential conjecture that has initiated and contributed to several research efforts in the field is Viterbo's volume-capacity conjecture \cite{viterbo2000}:
 \begin{conjecture}[Viterbo, 2000]
\label{viterbo_conj}
For any symplectic capacity $c$ and any convex domain $K \subset \R^{2n}$,
$$ c(K)^n \leq n! \mathrm{Vol}(K).$$
\end{conjecture}
Note that since symplectomorphisms are volume preserving, 
Conjecture \ref{viterbo_conj} is trivially satisfied for the Gromov width.
A stronger version of this conjecture asserts that all symplectic capacities coincide on the class of convex domains. 

\medskip 
Since its introduction, Viterbo's conjecture has inspired numerous related works, many of which provide substantial evidence supporting its validity.
In particular, one can verify that whenever the supremum over ball embeddings in $c_{\mathrm{G}}(K)$ is attained for a convex domain $K$, one has  $c_{\mathrm{G}}(K) = \ehzcap(K),$ thereby confirming Conjecture~\ref{viterbo_conj} for $\ehzcap$ in this case. 
The argument relies on the observation that if the intersection $\partial B \cap \partial K$ (where $B$ is the embedded ball) contains an entire closed characteristic of $B$, 
 then this curve is also a closed characteristic of 
$K$, and hence
the two capacities coincide. If, on the contrary, no closed characteristic of the ball lies entirely on $\partial K$, one can, roughly speaking, push the embedded ball further into the interior via a Hamiltonian diffeomorphism and thus embed a larger ball (cf. \cite[Theorem 1.2]{fathi22} and references within). 
Another folklore fact is that the strong Viterbo's conjecture, i.e. that all capacities coincide, holds for convex domains with $S^1$-symmetry (cf.~\cite{gutt-hutchings-ramos}).
Moreover, while for certain capacities the capacity-volume ratio is unbounded on some classes of non-convex domains, such as star-shaped ones \cite{hermann}, it was shown in \cite{ArtsteinAvidanOstroverMilman} that the conjecture holds up to a universal constant.
In \cite{abhs} (cf.~\cite{abbondandoloBenedetti}) it was proven that Conjecture~\ref{viterbo_conj} holds
in a $C^3$-neighborhood of the 
ball, and moreover that every smooth local maximizer in $\R^4$ which is strictly convex is symplectomorphic to the ball.  Later, in \cite{edtmair} (cf.~\cite{abbondandolo-benedetti-edtmair}) 
it was shown that the aforementioned stronger version of the conjecture is valid within a 
$C^2$-neighborhood of the ball.
In addition, Viterbo's conjecture was shown to hold for specific examples.
The papers \cite{balitsky, rudolf} prove Conjecture \ref{viterbo_conj} for certain Lagrangian products.
The conjecture was also verified for  monotone toric domains in~\cite{gutt-hutchings-ramos}, and for specific domains associated with Hamiltonians of classical mechanical type in 
\cite{karasev-sharipova}.
In \cite{capacity_mahler} it is shown that the famous Mahler conjecture on the volume product of centrally symmetric convex bodies, which had remained open for the past 85 years, and was recently proven in  dimension three \cite{mahler3dim}, is equivalent to the restriction of Conjecture~\ref{viterbo_conj} to the Ekeland--Hofer--Zehnder capacity of Lagrangian products of centrally symmetric convex bodies $K \subset \R^n_q$ and their dual bodies  $K^\circ \subset \R^n_p$.
We also observe that employing algorithms from \cite{pazit} and \cite{Ch-H} to compute the Ekeland--Hofer--Zehnder capacity led to the identification of several new examples with the same capacity-volume ratio as the ball.
Our main result is:

\begin{theorem}
\label{counterexample_thm}
Viterbo's conjecture fails 
for every 
$n \geq 2$.
\end{theorem}

The proof relies on the close connection between the Ekeland--Hofer--Zehnder capacity of Lagrangian products of convex domains in the classical phase space and Minkowski billiard dynamics~\cite{AA-O}. Minkowski billiards, introduced by Gutkin and Tabachnikov in~\cite{GutTab}, provide a natural generalization of classical billiards to the Finsler setting, where the Euclidean structure is replaced by a Minkowski distance function. In this framework, we consider two convex bodies \( K, T \subset \mathbb{R}^n \), which for simplicity we assume to contain the origin.
 One of these bodies, say \( K \), serves as the billiard table, while the other body \( T \) induces a potentially asymmetric norm via the support function  
$
h_T(u) := \sup \{ \langle x, u \rangle \mid x \in T\}
$
which governs the billiard dynamics within \( K \). Specifically, the reflection law for Minkowski billiards arises from a variational principle applied to the length functional, where the Euclidean length is replaced by \( h_T \). In what follows, we refer to the resulting billiard trajectories as \( T \)-billiard trajectories in \( K \).
Furthermore, since the function \( h_T \) is the gauge function of the dual body \( T^{\circ} = \{ x \in {\mathbb R}^{n} \, | \, \langle x,y \rangle \leq 1, \ \forall y \in T\} \), we will refer to the corresponding distance associated with the support function \( h_T \) as the \( T^{\circ} \)-length. 
In particular, for a  closed $T$-billiard trajectory in $K$ described by a  polygonal curve $\gamma$ with vertices $q_1,q_2,\ldots,q_m$, where $q_i \in \partial K$, the $T^\circ$-length of $\gamma$ is given by
$$ {\rm Length}_{T^{\circ}} (\gamma) = \sum_{i=1}^{m} h_T (q_{i+1}-q_i),$$
where $q_{m+1} = q_1$. As will be explained in more detail in Section~\ref{Section-Mink-Bill} below, the Ekeland--Hofer--Zehnder capacity of the Lagrangian product of two convex domains \( K \times T \) is given by the \( T^{\circ} \)-length of the shortest periodic \( T \)-billiard trajectory in \( K \) (cf. Section 2.4 in~\cite{AA-O}, and Section 1.2 in~\cite{Mink-Bill-Rudolf}).
\medskip

\begin{figure}[H]
\centering
\includegraphics[width=\linewidth]{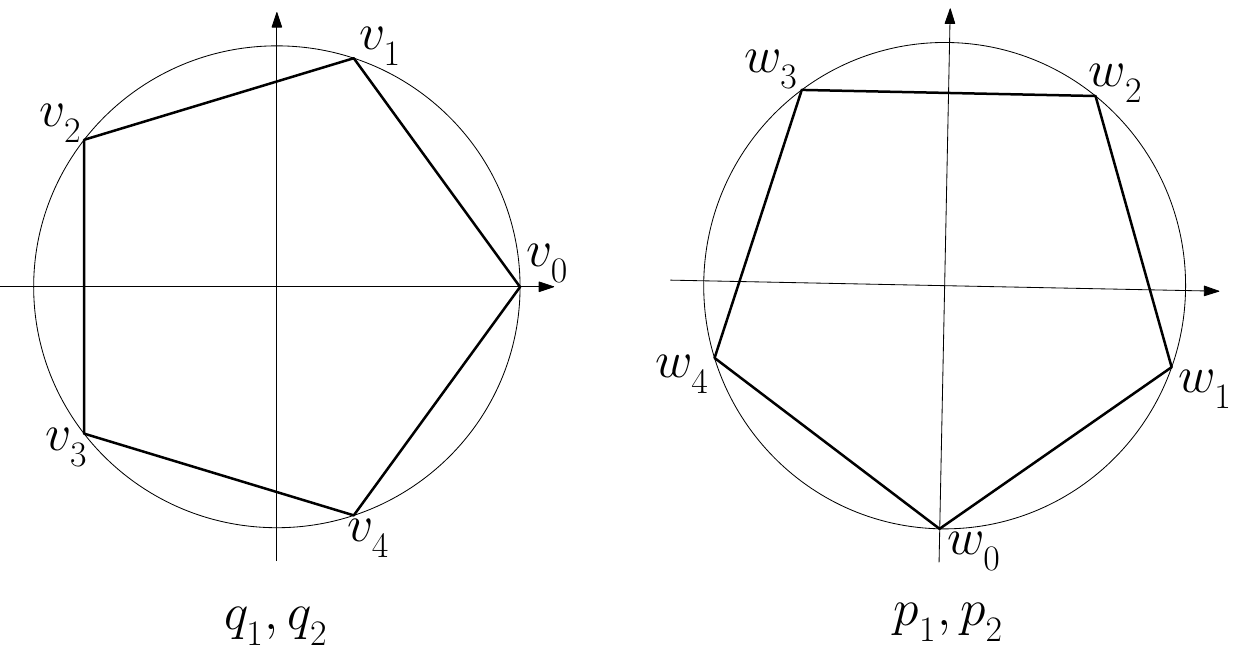}
\caption{The pentagon $K$ times the rotated pentagon $T$.}
\label{pentagonFig}
\end{figure}

Let us present an explicit counterexample to Viterbo's conjecture.
Consider the Lagrangian splitting of $\R^4 = \R^2_q \oplus \R^2_p$. Let $K \subset \R^2_q$ be the regular pentagon with vertices $\{(\cos(\frac{2\pi k}{5}), \sin(\frac{2\pi k}{5}))\}_{k=0}^4$,
and let $T \subset \R^2_p$ be a rotation of $K$ by $90^o$, i.e., the vertices of $T$ are $\{(\cos(-\frac{\pi}{2} + \frac{2\pi k}{5}), \sin(-\frac{\pi}{2} + \frac{2\pi k}{5}))\}_{k=0}^4$ (see Figure \ref{pentagonFig}). The main ingredient in the proof of Theorem~\ref{counterexample_thm} is the following proposition.

\begin{proposition}
\label{counterexample_prop}
The Ekeland--Hofer--Zehnder capacity of the Lagrangian product $K \times T$ is the $T^{\circ}$-length of any one of the $2$-bounce $T$-billiard trajectories along a diagonal of the pentagon $K$, i.e.,
$$ \ehzcap(K \times T) = 2 \cos(\frac{\pi}{10}) \Bigl (1 + \cos(\frac{\pi}{5}) \Bigr).$$
\end{proposition}

We observe that the  Reeb dynamics of \( \partial (K \times T) \) exhibit interesting global behavior. In particular, there is an open and dense set such that the Reeb trajectories passing through it are all closed. This phenomenon is examined in more detail in Section~\ref{sec-reeb-dyanmics} below.

\medskip

In what follows, for a convex domain \( A \subset {\mathbb R}^{2n} \), we denote   
the symplectic systolic ratio of \( A \) by
\[
{\mathrm {Sys}}(A) := \frac {\ehzcap(A)^n}{n! {\rm Vol}(A)}. \]
Using this notation, Viterbo's conjecture for $\ehzcap$ can be stated as \( \mathrm{Sys}(A) \leq 1 \) for all convex domains \( A \subset {\mathbb R}^{2n} \). Furthermore, we shall also make use of the notion of \( p \)-product. Recall that for two subsets \( X \subset {\mathbb R}^m \) and \( Y \subset {\mathbb R}^n \) containing the origin, and for \( 1 \leq p \leq \infty \), their \( p \)-product is defined as

\[
X \otimes_p Y := \bigcup_{0 \leq t \leq 1} \left( (1-t)^{1/p} X \times t^{1/p} Y \right) \subset \mathbb{R}^{m} \times \mathbb{R}^{n}.
\]  

Note that in the limit case \( p=\infty \), this operation corresponds to the usual product, whereas for \( p=1 \), it results in the convex hull \( {\rm Conv}\{( K \times \vec{0})  \cup (\vec{0} \times T)\} \).

\begin{proof}[{\bf Proof of Theorem \ref{counterexample_thm}}]
For the case \( n=2 \), since the area of both \( K \) and \( T \) is given by  
$
A = \frac{5}{2} \sin\left(\frac{2\pi}{5}\right)
$, Proposition \ref{counterexample_prop} immediately implies   
\[
\mathrm{Sys}(K \times T) = \frac{\ehzcap(K\times T)^2}{2 A^2} = \frac{\sqrt{5} + 3}{5} > 1.
\]  
To extend this counterexample to arbitrary dimensions, we utilize the fact (see Theorem 1.2 in \cite{p-prod}) that for any convex domains $X \subset \R^{2m}$ and $Y \subset \R^{2n}$ satisfying \( \ehzcap(X) = \ehzcap(Y) \), one has that $\mathrm{Sys}(X \otimes_2 Y) = \mathrm{Sys}(X) \cdot \mathrm{Sys}(Y).$
The proof of Theorem \ref{counterexample_thm} is then completed by taking the symplectic 2-product of \( K \times T \) with a Euclidean ball $B^{2n-4}(\sqrt{\ehzcap(K \times T)/\pi})$ of equal capacity.
\end{proof}

\begin{remark}
Using the formula provided in Theorem 1.1 of~\cite{pazit}, we have numerically computed  $\ehzcap( K \times T)$ using a MATLAB  program on standard personal computers, with a running time of less than one minute. 
The numerical computation yielded the claimed capacity value for the product $K \times T$ with a negligible numerical error. An  implementation of the algorithm, which reproduces our results, is available in~\cite{pazit-website}. Alternative methods for conducting the above numerical computation are discussed in~\cite{Ch-H} and~\cite{R-K}. For completeness, we  present an independent proof of Proposition~\ref{counterexample_prop} below.
\end{remark}

\begin{remark}
A notable class of domains in \(\mathbb{R}^{2n}\), known as toric domains, consists of regions that are invariant under the standard \(n\)-torus action on \(\mathbb{C}^n\). These domains are of special interest in symplectic geometry, as their symmetric features facilitate the analysis of various symplectic and dynamical properties (see, e.g.,~\cite{Gu-Ha,gutt-hutchings-ramos,C-CG-F-H-R,CG} and the references therein). 
A natural question is, roughly speaking, how far a general convex domain is from the class of toric domains, and more specifically, whether every convex domain admits a symplectic image that is a star-shaped toric domain.  
In a recent paper, Hutchings~\cite{Hut24} established that the answer to this question is negative in \(\mathbb{R}^4\). As J. Gutt pointed out to us, Theorem~\ref{counterexample_thm} offers an alternative proof of this result.
Indeed, in~\cite{gutt-hutchings-ramos}, it was shown that all (normalized) symplectic capacities coincide for a subclass of toric domains known as \textit{monotone toric domains} in \(\mathbb{R}^4\). This result was later extended to higher dimensions in~\cite{CG-H}.  
Consequently, since any star-shaped toric image of a convex domain in \(\mathbb{R}^4\) must be monotone~\cite[Proposition 1.8]{gutt-hutchings-ramos}, it follows that the Lagrangian product \(K \times T\),  after a smoothing procedure, cannot be symplectomorphic to any star-shaped toric domain.
\end{remark}

\subsection{ Discussion and Open Questions} 

While Theorem~\ref{counterexample_thm} provides a counterexample to Viterbo's volume-capacity conjecture, it also gives rise to several new questions. In this subsection we explore some 
 of these questions.

\medskip

\noindent {\bf (i) Viterbo's Conjecture for  Symmetric Convex Domains}

\vskip 3pt 

As mentioned earlier, substantial evidence suggests the possibility of an underlying fundamental reason that accounts for the validity of Viterbo's conjecture in many special cases.
While Theorem \ref{counterexample_thm} demonstrates that Conjecture \ref{viterbo_conj} does not hold universally for the class of convex domains,
it remains an intriguing question to determine
the appropriate subclass for which Conjecture \ref{viterbo_conj} is applicable. Specifically, 

\vskip 3pt
 
 {\bf Question 1:} Does Viterbo's conjecture hold for the class of centrally symmetric convex domains?

\vskip 3pt

Viterbo's conjecture is known to hold for convex domains that exhibit sufficient symmetry (see, e.g., 
Proposition 1.4 in~\cite{gutt-hutchings-ramos}). 
This naturally leads to the question of how much the symmetry assumptions can be relaxed while still ensuring the conjecture's validity. 
In particular, the subclass of centrally symmetric convex domains is of special interest, as a positive resolution to the aforementioned question would  imply the validity of Mahler's conjecture \cite{capacity_mahler}.

\medskip 

\noindent {\bf (ii) The Maximal Systolic Ratio}

\vskip 3pt

Another natural question that arises from Theorem~\ref{counterexample_thm} is to determine the exact maximal systolic ratio among 
convex domains. More precisely, let $\mathcal K^{2n}$ denote the class of convex domains in $\R^{2n}$, and denote the 
maximal symplectic systolic constant of this class by  $ \mathrm{Sys}(\mathcal K^{2n}) := \sup_{K \in \mathcal K^{2n}} \mathrm{Sys}(K). $

\medskip 
 
{\bf Question 2:} Determine the value $\mathrm{Sys}(\mathcal K^{2n})$.

\medskip

In contrast with the larger class of star-shaped domains in which the capacity-volume ratio is unbounded~\cite{hermann}, for convex domains one has $\mathrm{Sys}(\mathcal K^{2n}) \leq \alpha^n$, for some implicit dimension-independent universal constant $\alpha > 0$ (see \cite{ArtsteinAvidanOstroverMilman}). 
Moreover, it follows from~\cite{p-prod} that 
$\mathrm{Sys}(\mathcal K^{2n}) \leq \mathrm{Sys}(\mathcal K^{2m})$ whenever $n \leq m$.
Using John ellipsoids and compactness with respect to the Hausdorff metric, one can check that the supremum 
 in the definition of $\mathrm{Sys}(\mathcal K^{2n})$ 
 is attained, and can be replaced by a maximum.
Theorem \ref{counterexample_thm} naturally raises the question of identifying the maximizers, and exploring their dynamical properties. 

\medskip

We note that the largest systolic ratio currently known in ${\mathbb R}^{2n}$ is $ \Bigl( \tfrac{\sqrt{5} + 3}{5} \Bigr)^{\lfloor \frac{n}{2} \rfloor}.$ Indeed, this value can be achieved in a similar manner to the proof of Theorem~\ref{counterexample_thm}, which we now describe. 
When \( n \) is even, the construction is done by taking the symplectic 2-product of \( \frac{n}{2} \) copies of \( K \times T \). If \( n \) is odd, it instead involves the 2-product of \( \frac{n-1}{2} \) copies of \( K \times T \) along with a 2-dimensional ball of the same capacity. In both cases, \( K \) and \( T \) refer to the pentagons described in Proposition~\ref{counterexample_prop}.
Since the 2-product of two polytopes is not a polytope, the above construction yields a counterexample that is not a polytope. Nevertheless, we remark that in every dimension \( 2n \), there exists a Lagrangian product of convex polytopes that provides a counterexample to Viterbo's conjecture  with the aforementioned systolic ratio. 
Indeed, consider the following domains
\begin{align*}
K_{2k} = \underbrace{K \otimes_1 \cdots \otimes_1 K}_{k \text{ times}},& 
\quad K_{2k+1} = K_{2k} \otimes_1 [-l,l],\\
\quad T_{2k} = \underbrace{T \otimes_\infty \cdots \otimes_\infty T}_{k \text{ times}},&
\quad T_{2k+1} = T_{2k} \otimes_\infty [-l,l],
\end{align*}
where $l = \sqrt{\ehzcap(K \times T)}/2$.
Now ${\rm Sys}(K_n \times T_n) = \Bigl( \frac{\sqrt{5} + 3}{5} \Bigr)^{\lfloor \frac{n}{2} \rfloor}$ follows from the fact~\cite[Proposition 5.1.4]{pazit-thesis} that for convex domains $X_1,Y_1 \subset {\mathbb R}^{m}$, and $X_2,Y_2 \subset {\mathbb R}^{n}$ containing the origin with $\ehzcap(X_1 \times Y_1) = \ehzcap(X_2 \times Y_2)$, one has
$$ {\rm Sys} \bigl ((X_1 \otimes_1 X_2) \times (Y_1 \otimes_{\infty} Y_2) \bigr ) = {\rm Sys} (X_1 \times Y_1) \cdot {\rm Sys}(X_2 \times Y_2).$$

\begin{remark} Another notable class of domains in $\R^{2n}$, which includes convex domains, is the class $\mathcal {DC}^{2n}$ of dynamically convex domains (see~\cite{HWZ-3dconvex}).
In response to a longstanding open question, it was recently proved (first in $\R^4$~\cite{dynconvnotconvdim4}, and then in every dimension~\cite{dynconvnotconvdim2n}) that 
there exist dynamically convex domains in $\R^{2n}$
that are not symplectomorphic to any convex domain.
The maximal systolic constant for this class $\mathrm{Sys}(\mathcal {DC}^{2n})$ can be defined in a similar manner, where one replaces the Ekeland--Hofer--Zehnder capacity with the minimal action among closed characteristics on the boundary. In~\cite{Systolicratio} it was shown that $\mathrm{Sys}(\mathcal {DC}^{4}) \geq 2$. To the best of our knowledge, it is currently unknown if this quantity is bounded from above.  
It is a natural question  to determine how different $\mathrm{Sys}(\mathcal {DC}^{2n})$ and $\mathrm{Sys}(\mathcal {K}^{2n})$ are.
\end{remark}

\noindent {\bf (iii) Different Capacities of Convex Domains }

\vskip 3pt

As previously mentioned, a well-known stronger version of Viterbo's volume-capacity conjecture suggests that all (normalized) symplectic capacities coincide on convex domains, i.e., \( c_{\mathrm{G}}(K) = c_{\mathrm{Z}}(K) \) for all convex \( K \in {\mathcal K}^{2n} \) (see, e.g.,~\cite{capacity_survey_1, hermann}). 
Although Theorem~\ref{counterexample_thm} immediately disproves this conjecture, the fact that many symplectic capacities do coincide on this class (see~\cite{A-E-K,hofer-zehnder, hofer-ekeland, viterbo_cap, irie, abbondandolo-kang, Gu-Ha, ginzburg-shon, gutt-hutchings-ramos}), naturally leads to the following question:

\medskip

{\bf Question 3:} How many distinct\footnote{By distinct, we mean normalized capacities which cannot be represented as a homogeneous and monotone function of the other normalized capacities, as discussed in Section~3.3 of~\cite{capacity_survey_1}.} (normalized) symplectic capacities exist on the class of convex domains in ${\mathbb R}^{2n}$?

\medskip

We remark that in dimension $4$, it seems that currently there are only two known distinct capacities: $c_{\mathrm{G}}$ and $c_{\mathrm{Z}}$. 
 Indeed, it was recently shown in~\cite{A-E-K} that the Ekeland--Hofer--Zehnder capacity coincides with the cylindrical capacity for all convex domains in \( \mathbb{R}^4 \). 
Furthermore, combining results from~\cite{edtmair} and~\cite{Hryniewicz-Hutchings-Ramos}, one concludes that the first embedded contact homology  capacity and the cylindrical capacity coincide for the larger class of dynamically convex domains in $\R^4$.

\medskip

We finish this section with the following  remark: 

\begin{remark}
An immediate corollary of Theorem~\ref{counterexample_thm} is the existence of convex domains with a systolic ratio of one that are not symplectomorphic to the Euclidean ball, contrary to a conjecture in~\cite{viterbo2000}. Indeed, we owe the following argument to E. Kerman: From Theorem~\ref{counterexample_thm}  and a standard approximation argument it follows that for every $n \geq 2$ there is a smooth convex body $C \subset {\mathbb R}^{2n}$ with ${\rm Sys}(C) > 1$.
Let $L \subset {\mathbb R}^{2n}$ be a convex domain such that ${\rm Sys}(L) < 1$,
and consider $\lambda L + (1 - \lambda)C$. 
By continuity, for some $\lambda' \in (0,1)$ 
the domain $\lambda' L + (1-\lambda')C$ has systolic ratio one. 
Moreover, \(\lambda'\) can be chosen so that even an arbitrarily small increase in its value causes the systolic ratio of the domain to exceed one.  
From~\cite[Corollary 2]{abbondandoloBenedetti}, it follows that any symplectic image of the Euclidean ball admits a \(C^3\)-neighborhood in which every convex domain satisfies Viterbo’s conjecture.  
Consequently, the convex combination associated with \(\lambda'\) cannot be symplectomorphic to the Euclidean ball.
\end{remark}

\medskip

\noindent{\bf Acknowledgment:} The authors are grateful to Alberto Abbondandolo, Leonid Polterovich, and Shira Tanny for insightful feedback and suggestions. We also sincerely thank the referees for their extremely valuable remarks and comments, which have significantly contributed to improving the manuscript. Both authors were partially supported by the ISF grant No. 938/22. 
P.HK. was partially supported by a grant from the Institute for Advanced Study by the Erik Ellentuck Fellow Fund.

\section{Minkowski billiards}  \label{Section-Mink-Bill}

In this section we recall the definition of Minkowski billiards. 
For more information on Minkowski billiard dynamics see, e.g.,~\cite{AA-O,GutTab},
and~\cite{Mink-Bill-Rudolf}, where no smoothness assumptions are made on the convex bodies involved.

\medskip

\begin{figure}[H]
\centering
\includegraphics[width=\textwidth]{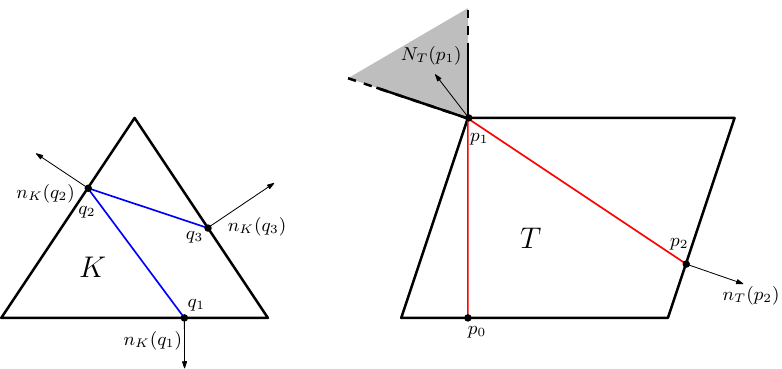}
	 	\caption{An example of a $T$-billiard trajectory in $K$.}
	 	\label{billiard-fig}
\end{figure}

The study of billiard dynamics in Minkowski geometries was initiated by Gutkin and Tabachnikov~\cite{GutTab} as an important special case of the natural extension of Euclidean billiard dynamics to the
Finsler setting. The connection between convex Minkowski billiards and the Ekeland–Hofer–Zehnder capacity was established in~\cite{AA-O} (cf.~\cite{Mink-Bill-Rudolf, R-K, viterbo2000}), and subsequently utilized in~\cite{capacity_mahler} to relate Viterbo's conjecture and Mahler's conjecture from convex geometry.  From the perspective of geometric optics, Minkowski billiard trajectories describe the propagation of light in a homogeneous anisotropic medium with perfectly reflecting boundaries. 
We focus on the specific case of Minkowski billiards within a convex body $K$.
Endow \( K \) with a (not necessarily symmetric) distance function 
associated with a convex body $T$, and consider billiard dynamics in \( K \) with respect to the induced geometry. 
More precisely, let \( K \subset \mathbb{R}^{n}_q \) be a convex domain (the billiard table), and let 
 \( T \subset \mathbb{R}^{n}_p \) be a convex domain which has the origin in its interior. The corresponding Minkowski geometry induced by $T$ is given by the support function 
$
h_T(u) := \sup \{ \langle x, u \rangle \mid x \in T\},
$ which induces a (possibly asymmetric) distance function 
$d_{T^\circ}(x,y) = h_T(x-y)$. 
This notation is motivated by the fact that the polar body $T^\circ$ is the unit ball with respect to this distance function.
Note that when $T$ is the Euclidean unit ball, $d_{T^\circ}$ is the usual Euclidean metric. When \( K \) and \( T \) are strictly convex smooth domains, the corresponding billiard dynamics in \( K \) is determined by the reflection associated with the distance function \( d_{T^\circ} \). That is, for fixed points \( x, z \in \partial K \), the points \( x, y, z \), where \( y \in \partial K \), form a \( T \)-billiard trajectory if \( y \) is an extremal point of the function \( d_{T^\circ}(x,y) + d_{T^\circ}(y,z) \). Since we do not wish to assume smoothness or strict convexity, we provide the following definition (see Figure~\ref{billiard-fig}). Recall that for a convex domain \( K \subset \mathbb{R}^{n} \), 
the outer normal cone to \( K \) at \( x \in \partial K \) is defined as
\[
N_K(x) = \left\{ u \in \mathbb{R}^n \mid \langle u, y - x \rangle \leq 0 \text{ for all } y \in K \right\}. 
\]
We note that when $\partial K$ is smooth $N_K(x) = \R_+ n_K(x)$, where $n_K(x)$ is the unique unit outer normal.

\begin{definition}[Minkowski Billiards] \label{def-T-billiard-in-K}
  A closed polygonal curve with vertices \( q_1, \dots, q_m \), \( m \in \mathbb{N}_{\geq 2} \), on \( \partial K \) is a closed $T$-billiard trajectory in $K$ if
there are points \( p_1, \dots, p_m \) on \( \partial T \) such that for every $ 1 \leq j \leq m$,
\[
\begin{cases}
q_{j+1} - q_j \in N_T(p_j), \\
p_{j+1} - p_j \in -N_K(q_{j+1}).
\end{cases}
\]
 \end{definition}

Here, and in what follows, 
indices \( 1, \dots, m \) are considered modulo \( m \) in \( \mathbb{Z} \).
Note that this definition naturally extends to polygonal curves that are not closed.
A direct computation shows that when \( T \) is a Euclidean ball, and $K$ is strictly convex with piecewise smooth boundary, the corresponding billiard dynamics follows the well-known law of geometric optics: the angle of incidence equals the angle of reflection.

 \medskip

In the Euclidean case, Bezdek and Bezdek~\cite{bezdek-bezdek} introduced a useful characterization of length-minimizing billiard trajectories, stating that closed billiard trajectories in $
K$ of minimal length correspond to length-minimizing closed polygonal curves that cannot be translated into the interior of 
$K$.
This statement, proved in Theorem 1.1, and in  particular Lemma 2.4, of~\cite{bezdek-bezdek}, was later extended to Minkowski billiards (see e.g.,~\cite[Theorem 2.1]{non_symmetric_mahler}, and~\cite[Theorem 1]{Mink-Bill-Rudolf}).
More precisely, denote by $\mathcal{P}$ the set of all polygonal curves which cannot be translated into the interior of $K$:
$$ \mathcal{P}(K) = \{(q_1,\ldots,q_m) \mid \{q_1,\ldots,q_m\} \not\subset {\rm int} (K)+t, \; \forall t \in \R^n\}.$$
The result demonstrates  
that a polygonal curve $(q_1,\ldots,q_m)$ is a $T$-billiard trajectory in $K$ with minimal $T^\circ$-length if and only if it is a curve in $\mathcal{P}(K)$ and it has the minimal $T^\circ$-length among all such curves.
In addition, the proofs ibid confirm that the minimum-length polygonal curve in $ P(K)$  can always be chosen to have at most $n+1$ vertices.
This characterization of the minimal periodic billiard trajectory is a key component in the proof of Theorem~\ref{counterexample_thm} below.

\medskip

To relate Minkowski billiards with the Ekeland--Hofer--Zehnder capacity, recall that a closed characteristic of a smooth convex domain $\Sigma \subset {\mathbb R}^{2n}$ is a loop $\gamma \subset \partial\Sigma$ such that $\dg(t)$ is a positive multiple of $J \nabla H(\gamma(t))$, where $H: \R^{2n} \to \R$ has $\Sigma$ as a sub-level set, and $J :\mathbb R^{2n} \rightarrow {\mathbb R}^{2n} $ is the standard complex structure $J^2 = -Id$.
Note that after rescaling, each such characteristic is an orbit of the Reeb vector field associated to the standard contact form on $\partial \Sigma$.
The action of a closed characteristic $\gamma$ is defined by ${\mathcal A}(\gamma) = \int_{\lambda} \gamma$, where $\lambda$ is a primitive of $\omega$. 
If $\Sigma$ is not smooth, one can define a generalized closed characteristic, roughly speaking, by changing the condition above to require that $\dg(t)$ is a positive multiple of a vector in $ J \partial H(\gamma(t))$, where $\partial H$ is the sub-gradient of $H$ (see, e.g.,~\cite{ Kunzle}).  Moreover, one can show  that the Ekeland--Hofer--Zehnder capacity of $\Sigma$ is given by 
 the minimal action over all generalized closed characteristics on the boundary $\partial \Sigma$. 
 For more precise details, see e.g.,~\cite{AA-O, Kunzle}. In the case where \(\Sigma\) is the Lagrangian product \(K \times T\) of two convex domains containing the origin, consider the Hamiltonian function
$$ H(q,p) = \max \{ g_K(q), g_T(p) \}, $$
where for a convex domain $\Sigma$ the  associated gauge function $g_{\Sigma}$ is defined by 
$ g_{\Sigma}(x) := \inf \{ r \mid x \in r \Sigma \}
$. 
The corresponding characteristic directions, defined by this Hamiltonian, at a point $(q_0,p_0) \in {\rm int}(K) \times \partial T$, belongs to $(N_T(p_0),\vec{0})$, which creates movements in a straight line until the trajectory hits a point $q_1$ on the boundary of $K$.
The movement induced by the Hamiltonian now continues in a similar manner in a straight line in $T$ along a direction belonging to $-N_K(q_1)$.
Furthermore, the action of a closed characteristic of \(K \times T\) can be naturally related to the billiard length measured with respect to \(h_T\). Specifically, denoting the vertices of a closed \(T\)-billiard trajectory in \(K\) by \((q_1, \ldots, q_m)\). This trajectory can be lifted to a generalized closed characteristic in ${\mathbb R}^n_q \oplus {\mathbb R}^n_p$ whose symplectic action equals the \(T^\circ\)-length of the billiard trajectory:
\[
\sum_{i=1}^{m} h_T(q_{i+1} - q_i).
\]
This relationship follows, roughly speaking, from the fact that the symplectic action of a loop \(\gamma\) can be expressed as the integral of the one-form \(\lambda = \sum_{i=1}^n p_i \, dq_i\) along  $\gamma$.

\medskip

Overall, one has that closed characteristics of $K \times T$ are in correspondence with $T$-billiards in $K$, which in turn are in correspondence with the minimal $T^\circ$-length polygonal curves which cannot be translated into the interior of $K$.
Moreover, the Ekeland--Hofer--Zehnder capacity equals the minimal action over all (generalized) closed characteristics and the action of a $K \times T$ closed characteristic equals the $T^\circ$-length of a $T$-billiard trajectory in $K$. This can be summarized as follows (see Theorem 2.13 in~\cite{AA-O},
and Theorem 1 in~\cite{Mink-Bill-Rudolf}). 

\begin{theorem} \label{thm-capacity-lagrangian-product}
For any two convex domains $K \subset {\mathbb R}^n_q$ and $T \subset {\mathbb R}^n_p$ one has
\begin{align*}
    \ehzcap(K \times T) &= \min\{\mathcal{A}(\gamma) \mid \gamma \text{ is a closed characteristic of } K \times T\}\\
    &= \min\{ {\rm Length}_{T^\circ}(\eta) \mid \eta\text{ is a closed } T\text{-billiard in } K\}\\
    &= \min\{ {\rm Length}_{T^\circ}(\eta) \mid \eta  \in \mathcal{P}(K)\text{ has at most } n+1\text{ vertices } \}.
\end{align*}
\end{theorem}

\section{The Reeb Dynamics of $K \times T$}

Building on the connection between Minkowski billiards and the Reeb flow associated with Lagrangian products of convex domains, as discussed in Section~\ref{Section-Mink-Bill}, we now proceed to examine the \( T \)-billiard trajectories in \( K \) for the pentagons product configuration described in the introduction.  
We begin by proving Proposition~\ref{counterexample_prop} by determining the minimal \( T^\circ \)-length of a billiard trajectory in \( K \), and continue with a discussion about the general behavior of $T$-billiard trajectories.

\begin{proof}[{\bf Proof of Proposition~\ref{counterexample_prop}}]
  
Denote the vertices of the regular pentagon inscribed in the unit disc $K \subset \R^2_q$ by $v_k := (\cos(\frac{2\pi k}{5}), \sin(\frac{2\pi k}{5}))$, and the vertices of the regular pentagon $T \subset \R^2_p$ by  $w_k:= (\cos(-\frac{\pi}{2} + \frac{2\pi k}{5}), \sin(-\frac{\pi}{2} + \frac{2\pi k}{5}))$, where $0 \leq k \leq 4$. 
As discussed in Section~\ref{Section-Mink-Bill}, and in particular Theorem~\ref{thm-capacity-lagrangian-product}, the computation of $\ehzcap(K \times T)$ can be carried out by finding the minimal $T^{\circ}$-length among periodic $T$-billiard trajectories in $K$. Furthermore, it suffices to consider only such trajectories with two or three bounce points which cannot be translated into the interior of $K$.
Recall that, by the definition of the support function, and the particular choice of $T$, the $T^\circ$-length of a vector $v \in \R^2_q$ is given by
\begin{equation} \label{eqn-T-dual-length} 
h_T(v) = \max_{0\leq k \leq 4} \langle v, w_k \rangle. 
\end{equation}
A direct computation shows that the $T^{\circ}$-length of any  of the $2$-bounce billiard trajectories along the diagonals of the pentagon $K$ is
$ 2 \cos(\frac{\pi}{10}) (1 + \cos(\frac{\pi}{5}))$. 
Any other $2$-bounce trajectory which cannot be translated into the interior is a trajectory connecting one of the vertices $v_k$ with a point 
$\widetilde q$ on the opposite edge, i.e., the edge connecting the vertices $v_{k+2}$ and $v_{k+3}$ counting indices modulo five.
We claim that the $T^\circ$-length of such $2$-bounce trajectory equals the $T^\circ$-length of any of the diagonal trajectories.
Indeed, it follows immediately from~\eqref{eqn-T-dual-length} that for a point $\widetilde q = \lambda v_{k+2} + (1 - \lambda)v_{k+3}$, where $\lambda \in [0,1]$, the $T^{\circ}$-length of the 2-bounce trajectory connecting $v_k$ and $\widetilde q$ is  
\begin{align*}
    h_T(v_k-\widetilde q)+h_T(\widetilde q - v_k) &= \langle v_k - \widetilde q, w_{k+1} \rangle + \langle \widetilde q - v_k, w_{k-1} \rangle.
\end{align*}
Note that since $ v_{k+3} - v_{k+2}$ is orthogonal to $w_{k+1} - w_{k-1}$, 
the above expression is independent of $\lambda$,
and one has that all such 2-bounce trajectories have the same $T^{\circ}$-length as the diagonal connecting $v_{k+3}$ and $v_k$.

\medskip 

Now consider $T^{\circ}$-length minimizing $3$-bounce trajectories.
Let $\gamma$ be 
a minimal $T^\circ$-length $3$-bounce trajectory, and denote its vertices by $x_1, x_2, x_3$.
Since this curve cannot be translated into the interior, two of its vertices lie on adjacent edges, and the third vertex lies on the opposite edge. 
Since $K$ and $T$ exhibit symmetry with respect to rotations with angle $\frac{2\pi}{5}$, one can assume without loss of generality that the vertices of $\gamma$ lie on the edges $[v_1, v_2]$, $[v_4 , v_0]$ and $[v_3, v_4]$ (see Figure~\ref{gamma1Fig}).

\medskip

Note that since $T$ is not centrally symmetric, one needs to consider both the clockwise and counterclockwise orientations of $\gamma$. 
Let us first address the clockwise orientation case.
\begin{figure}[h]
\centering
\includegraphics[width=0.5\textwidth]{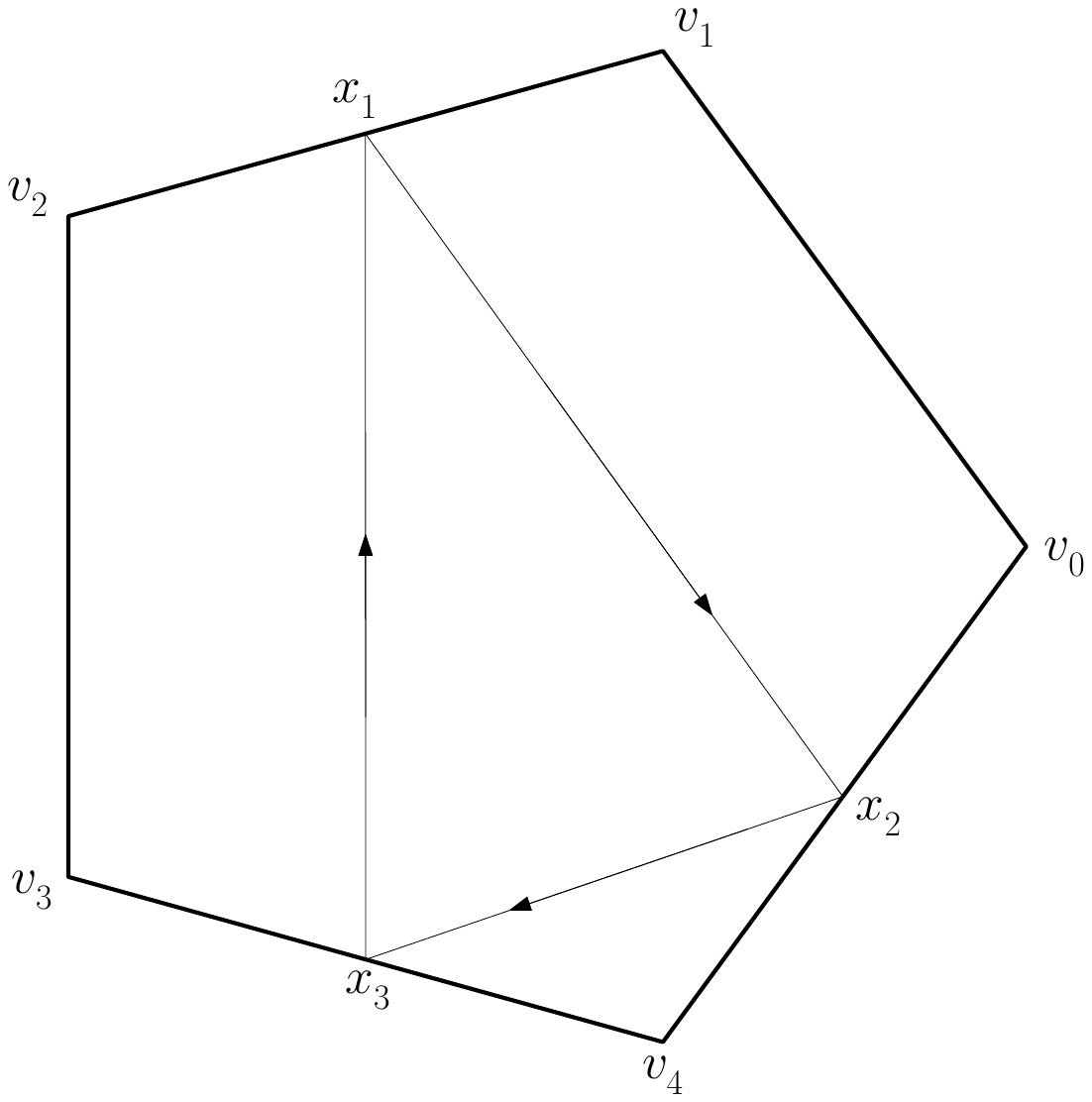}
	 	\caption{The billiard trajectory $\gamma$ in the clockwise orientation case.}
	 	\label{gamma1Fig}
\end{figure}
Let us fix $x_1$ and $x_3$ and  
explore the $T^\circ$-length of $\gamma$ when one moves $x_2$ on the line between $v_4$ to $v_0$.
Denote $x_2^\lambda := \lambda v_4 + (1-\lambda) v_0$, and let $\gamma_\lambda$ be the clockwise oriented polygonal curve with vertices $x_1,x_2^\lambda,x_3$.
It immediately follows from~\eqref{eqn-T-dual-length} that the $T^\circ$-length of the line from $x_2^{\lambda}$ to $x_3$ (for all possible choices of $x_2^\lambda$) is $h_T(x_3 - x_2^\lambda) = \langle x_3 - x_2^\lambda, w_4 \rangle$. 
For the line from $x_1$ to $x_2^\lambda$,  the value $h_T(x_2^\lambda - x_1)$ is either $\langle x_2^\lambda - x_1, w_1 \rangle$ or $\langle x_2^\lambda - x_1, w_0 \rangle$, depending on how $x_2^\lambda$ is chosen. 
Choose $\Lambda \in [0,1]$ such that the line from $x_1$ to $x_2^\Lambda$ is perpendicular to the line from $w_0$ to $w_1$.
We conclude that
\begin{align*}
h_T(\gamma_\lambda) = 
\begin{cases}
\langle x_2^\lambda - x_1, w_1\rangle + \langle x_3 - x_2^\lambda, w_4 \rangle + \langle x_1 - x_3, w_j \rangle, & \text{if } \lambda \leq \Lambda, \\
\langle x_2^\lambda - x_1, w_0\rangle + \langle x_3 - x_2^\lambda, w_4 \rangle + \langle x_1 - x_3, w_j \rangle, & \text{otherwise},
\end{cases}
\end{align*}
for some $0\leq j \leq4$.
Hence, the derivative of the $T^\circ$-length of $\gamma_\lambda$ with respect to $\lambda$ equals 
\begin{align*}
\frac{d}{d \lambda}h_T(\gamma_\lambda) = 
\begin{cases}
\langle v_4 - v_0, w_1 - w_4\rangle , & \text{if } \lambda \leq \Lambda, \\
\langle v_4 - v_0, w_0 - w_4\rangle , & \text{otherwise}.
\end{cases}
\end{align*}
Since $\langle v_4 - v_0, w_1 - w_4\rangle < 0$ and $\langle v_4 - v_0, w_0 - w_4\rangle = 0$, we get that the minimum is attained at $\lambda \in [\Lambda,1]$,
so we can assume that $x_2=v_4$. By the triangle inequality, the $T^\circ$-length of $\gamma_\lambda$ is greater or equal than the $T^\circ$-length of the 2-bounce trajectory between $x_1$ and $v_4$. 
This 2-bounce trajectory cannot be translated into the interior of $K$, 
a case that has already been discussed above.

\medskip 

Finally, we consider the case where the orientation is counterclockwise, and we denote the 3-bounce trajectory in this case as $\widetilde{\gamma}_\lambda$.
A similar argument yields
\begin{align*}
\frac{d}{d \lambda}h_T(\widetilde{\gamma}_\lambda) = 
\begin{cases}
\langle v_4 - v_0, w_2 - w_3\rangle , & \text{if } \lambda \leq \widetilde{\Lambda}, \\
\langle v_4 - v_0, w_1 - w_3\rangle , & \text{otherwise},
\end{cases}
\end{align*}
where $\widetilde{\Lambda}$ is chosen so that the line from $x_3$ to $x_2^{\widetilde{\Lambda}}$ is perpendicular to the line from $w_1$ to $w_2$.
Again, the minimum is attained at $\lambda \in [\widetilde{\Lambda},1]$, hence we may assume $x_2 = v_4$, and the rest follows in the same way as in the previous case. 
This completes the proof of the proposition. 
\end{proof}

\begin{figure}[H]
\centering
\begin{minipage}{.48\linewidth}
 \includegraphics[width=\linewidth]{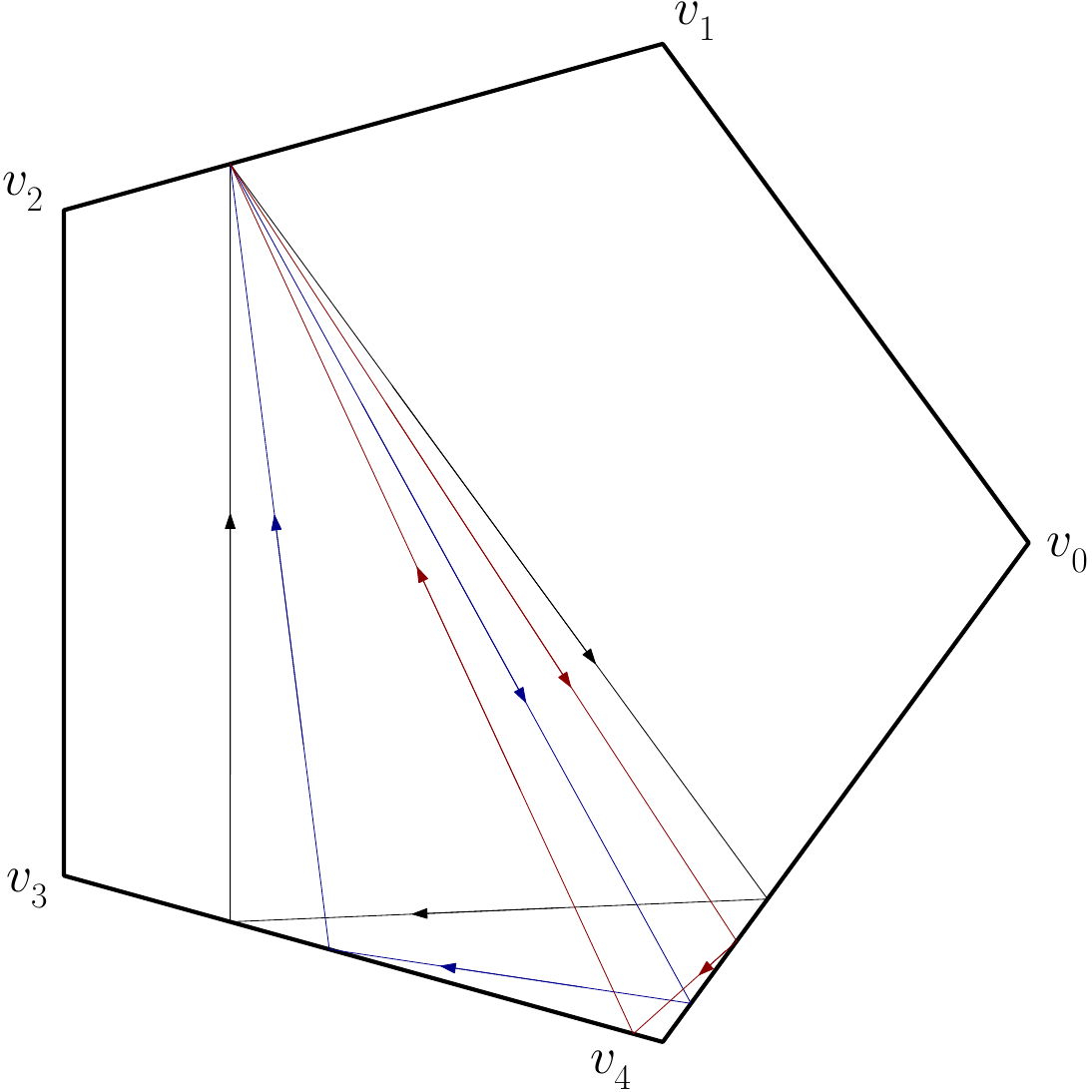}
\end{minipage}
\hspace{.01\linewidth}
\begin{minipage}{.48\linewidth}
 \includegraphics[width=\linewidth]{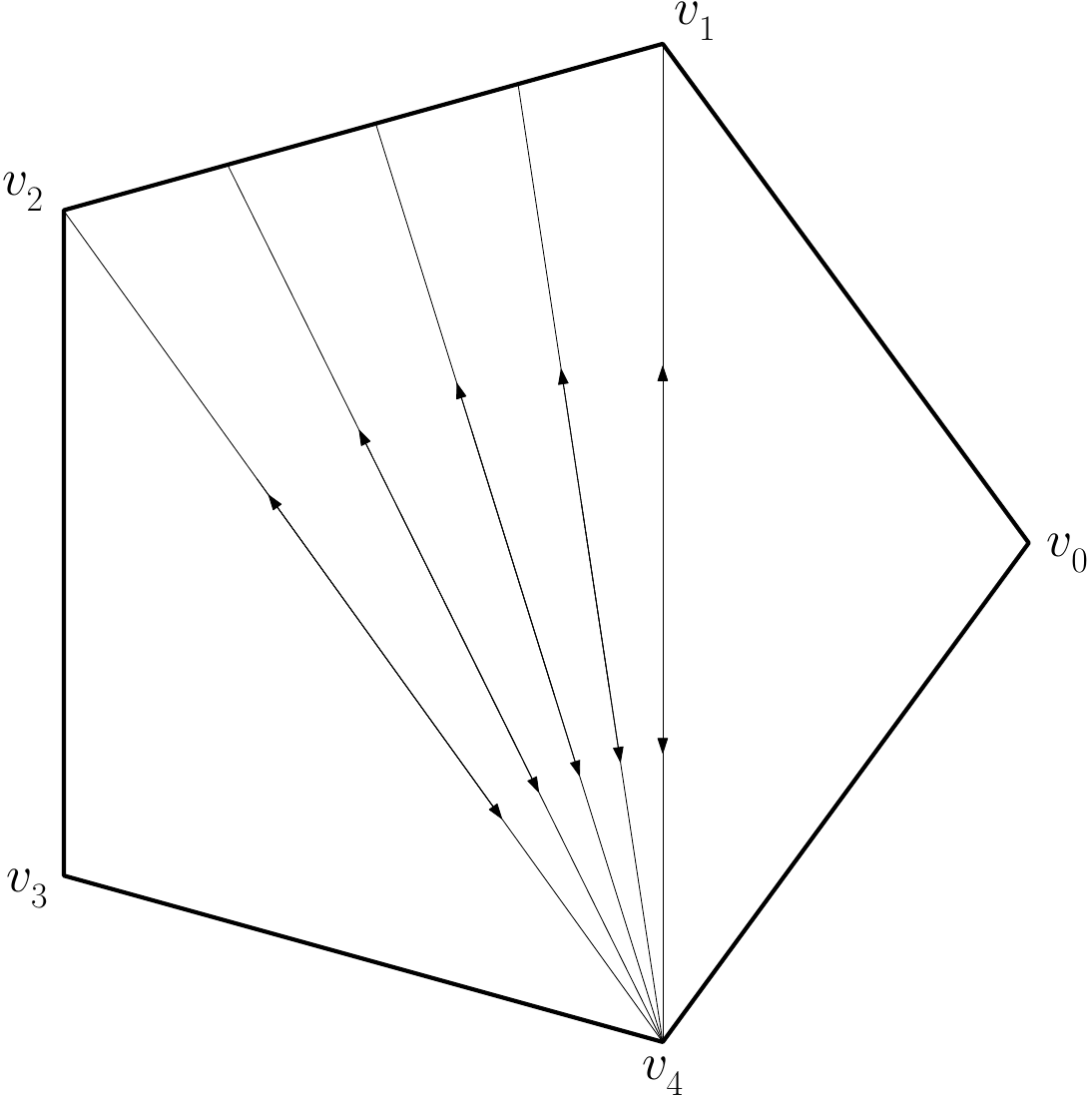}
\end{minipage}
\caption{Different minimal $T^\circ$-length billiard trajectories with three bounces (left) and two bounces (right).}
  \label{figure-minimizing-billiards}
\end{figure}

We remark that, in addition to the 2-bounce minimizing trajectories presented above, there also exist 3-bounce trajectories  with minimal \( T^\circ \)-length. Indeed, a variety of \( T^\circ \)-minimizing billiard trajectories can be obtained by following the arguments of the proof of Proposition~\ref{counterexample_prop}, and identifying equality cases where the triangle inequality has been used (which is possible due to the lack of smoothness of \( h_T \)).  
In particular, using the same notation as in the proof above, one can check that for a given \( x_1 \), all clockwise trajectories where \( x_2 \) and \( x_3 \) lie within the cone formed by the parallel lines to \( v_1v_0 \) and \( v_2v_3 \) passing through \( x_1 \) are billiard trajectories with minimal \( T^\circ \)-length (see  Figure~\ref{figure-minimizing-billiards}).

\subsection{The Generic Behavior}
\label{sec-reeb-dyanmics}
\begin{figure}[H]
\centering
\begin{minipage}{.7\linewidth}
 \includegraphics[width=\linewidth]{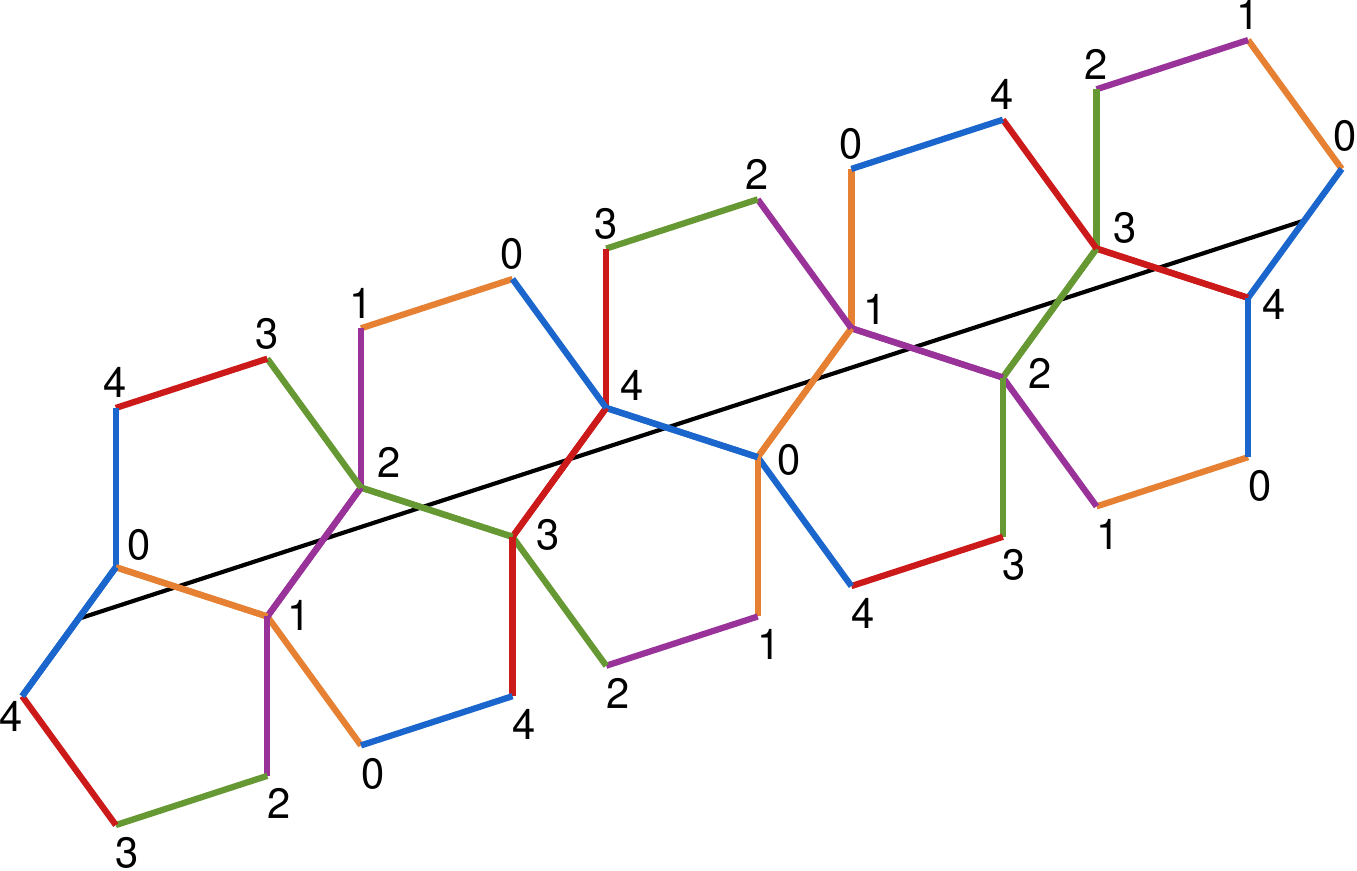}
\end{minipage}
\hspace{.01\linewidth}
\begin{minipage}{.25\linewidth}
 \includegraphics[width=\linewidth]{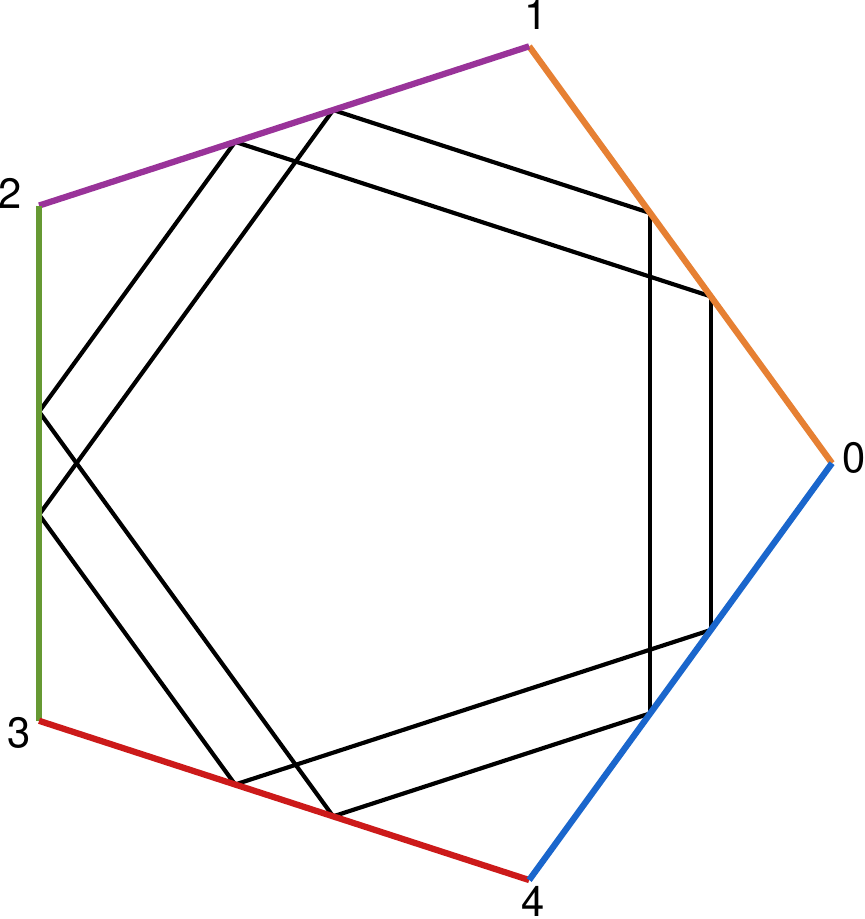}
\end{minipage}
\caption{A \( T \)-billiard trajectory in \( K \) of the first kind (right) and its associated unfolding (left).
}
  \label{figure_generic_billiard1}
\end{figure}

\begin{figure}[H]
\centering
\begin{minipage}{.75\linewidth}
 \includegraphics[width=\linewidth]{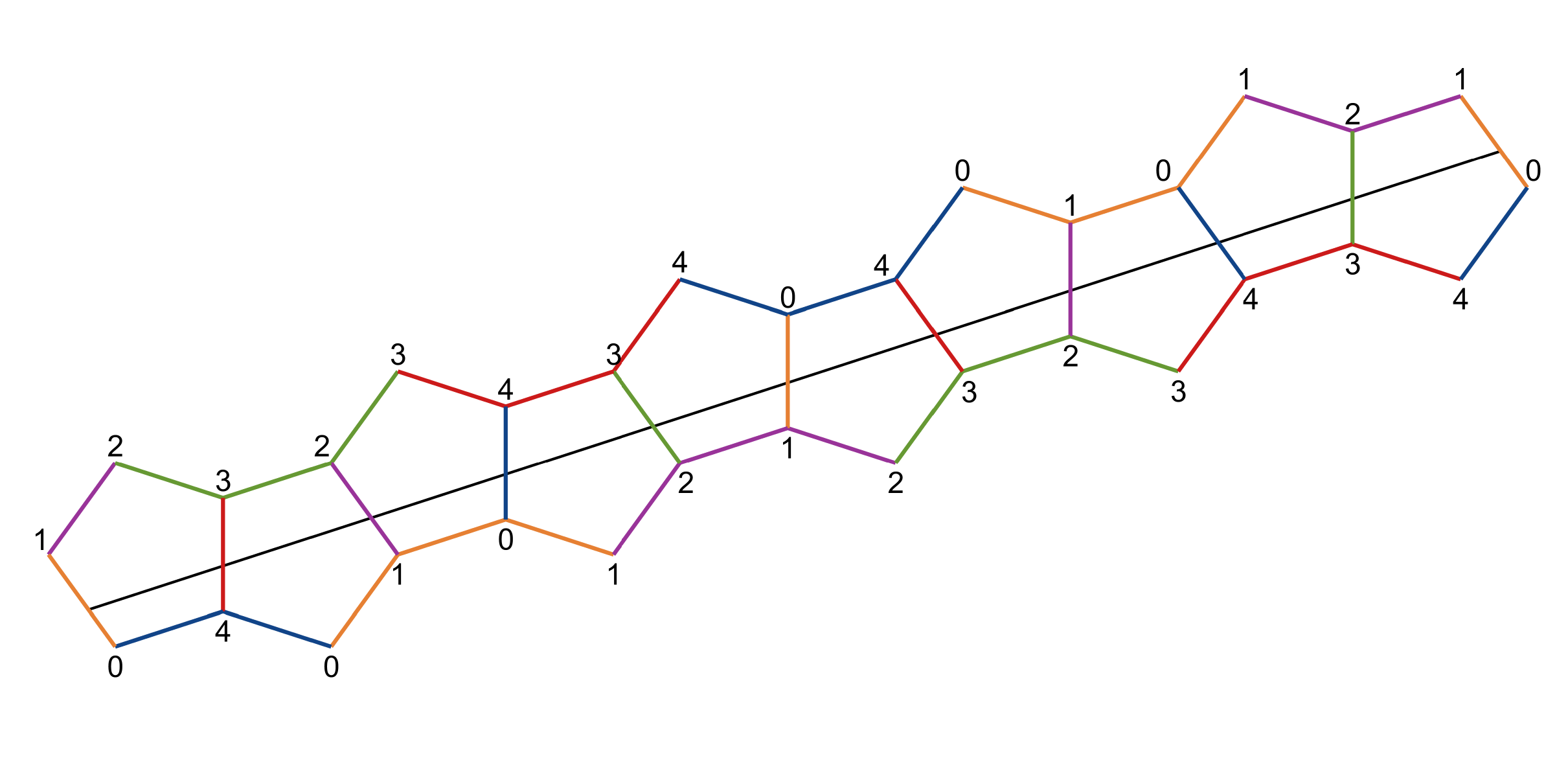}
\end{minipage}
\hspace{.01\linewidth}
\begin{minipage}{.22\linewidth}
 \includegraphics[width=\linewidth]{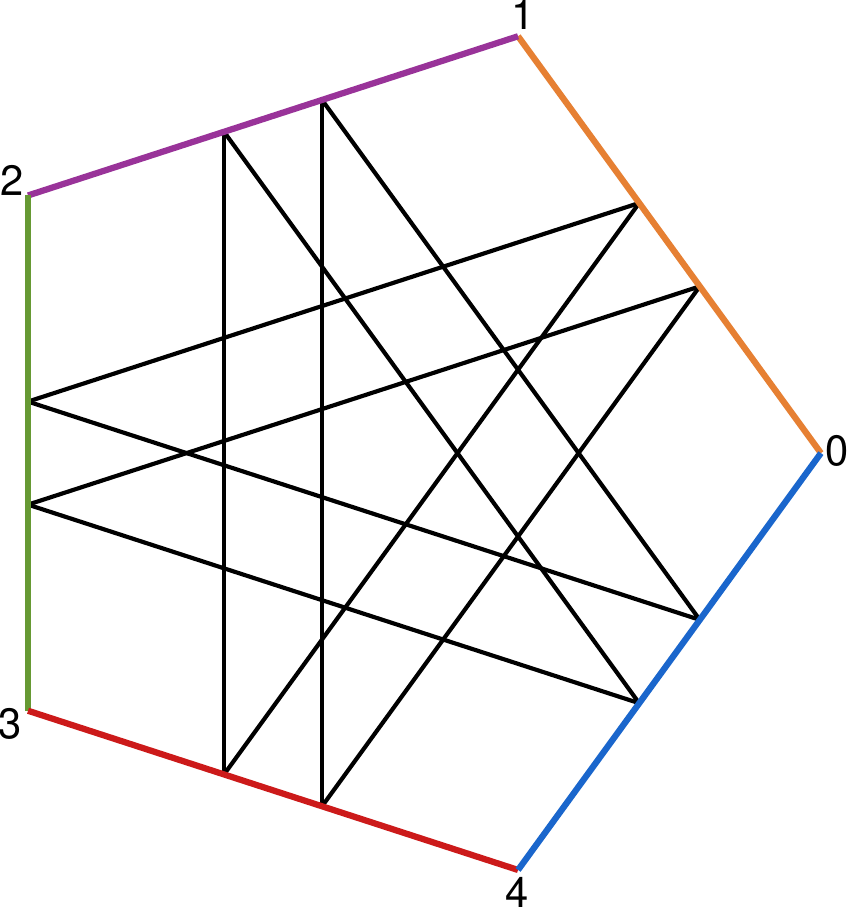}
\end{minipage}
\caption{A \( T \)-billiard trajectory in \( K \) of the second kind (right) and its associated unfolding (left).}
 \label{figure_generic_billiard2}
\end{figure}

In this subsection we consider the Reeb flow through a
point on the boundary \( \partial (K \times T) \),
corresponding to a \( T \)-billiard trajectory in \( K \) that avoids corners  and whose velocity vectors are always orthogonal to the edges of \( T \).
This type of behavior is generic in the sense that a small perturbation of the point ensures that it holds.
It turns out that these trajectories are periodic, and can be classified into two distinct types, each corresponding to one of two possible \( T^\circ \)-lengths.
Moreover, because the support functions of 
$T$ and its inscribed disc 
coincide in the five relevant directions, it follows that the reflection law governing these trajectories is the classical Euclidean one.
 Figure \ref{figure_generic_billiard1} and Figure \ref{figure_generic_billiard2} illustrate the two types of trajectories along with their unfolding, which demonstrate that their length is preserved under small perturbations.
 For each of the two types, there are also (non-generic in the sense described above) billiard trajectories of half the $T^\circ$-length passing through the middles of the edges, which close after 5 bounces instead of 10.
 We note that even these shorter trajectories are strictly longer than the diagonal trajectory described in Proposition~\ref{counterexample_prop}.

The dynamics in which periodic billiard trajectories run parallel to the edges through every point appear in the Lagrangian product of any regular polygon with its 90-degree rotation.
Such a behavior of the Reeb flow is closely related with the notion of ``symplectic Besse", i.e., 
a Reeb flow all whose orbits are closed.
We remark that any local maximizer of the capacity-volume ratios associated with the sequence of Ekeland--Hofer capacities, 
which is a \( 2n \)-dimensional smooth star-shaped domain, must have Besse boundary, and
moreover, in dimension $4$, such domains are symplectomorphic to rational ellipsoids (see~\cite{B-B-L-M}).
In the same paper, they pointed out non-smooth domains that have higher capacity-volume ratio than ellipsoids.
Further study of non-smooth domains with Besse-type dynamics would be particularly interesting.

\bibliography{references}
\bibliographystyle{siam}

\vskip10pt

\noindent Pazit Haim-Kislev \\
\noindent School of Mathematics, Institute for Advanced Study, Princeton, NJ, USA \\
\noindent e-mail: pazithaim@ias.edu
\vskip 10pt

\noindent Yaron Ostrover \\
\noindent School of Mathematical Sciences, Tel Aviv University, Israel \\
\noindent e-mail: ostrover@tauex.tau.ac.il

\end{document}